\documentclass[a4paper, 12pt, twoside]{article}
\usepackage{color}
\usepackage{amsmath}
\usepackage{amsthm}
\usepackage{amsfonts}
\usepackage{pifont}
\usepackage{bbm}
\usepackage{mathrsfs}
\usepackage{fancyhdr}
\usepackage{graphicx}

\usepackage{psfrag}
\usepackage{time}
\usepackage{txfonts}

\usepackage{stmaryrd}
\usepackage{mathrsfs}
\usepackage{txfonts}
\usepackage{amsfonts}

\usepackage{indentfirst}

\newcommand{\xiaosan}{\fontsize{15pt}{22pt}\selectfont}
\newcommand{\sihao}{\fontsize{14pt}{21pt}\selectfont}

\newcommand{\xiaosi}{\fontsize{12pt}{18pt}\selectfont}

\usepackage{curves}

\numberwithin{equation}{section}

\newtheorem{theorem}{ {Theorem}}[section]

\newtheorem{remark} {   {Remark}}[section]

\newtheorem{lemma} {  {Lemma}}[section]

\begin{document}
\setlength{\parindent}{2em}
\newpage
\fontsize{12}{22}\selectfont\thispagestyle{empty}
\renewcommand{\headrulewidth}{0pt}
 \lhead{}\chead{}\rhead{} \lfoot{}\cfoot{}\rfoot{}
\noindent

\title{\bf Absolutely Continuous Spectrum  for the Quasi-periodic Schr\"{o}dinger Operator in   Exponential Regime}
\author{{Wencai  Liu and Xiaoping Yuan*}\\
{\em\small School of Mathematical Sciences}\\
{\em\small  Fudan University}\\
{\em\small  Shanghai 200433, People's Republic of China}\\
{\small 12110180063@fudan.edu.cn}\\
{\small *Corresponding author: xpyuan@fudan.edu.cn}}
\date{}
\maketitle
\renewcommand{\baselinestretch}{1.2}
\large\normalsize
\begin{abstract}
Avila and  Jitomirskaya  prove that the quasi-periodic Schr\"{o}dinger operator  $H_{\lambda v,\alpha,\theta}$  has purely absolutely continuous spectrum for $\alpha $ in sub-exponential regime (i.e., $\beta(\alpha)=0$) with small $\lambda$, if $v$ is  real  analytic in a strip of real axis. In the present paper,  we   show that  for  all $\alpha$  with    $0<\beta(\alpha)<\infty$,  $H_{\lambda v,\alpha,\theta}$
has purely absolutely continuous spectrum with   small $\lambda$, if $v$ is real
    analytic in   strip  $  |\Im x|< C\beta $, where $C$ is a large absolute constant.
\end{abstract}

\setcounter{page}{1} \pagenumbering{arabic}\topskip -0.82in
\fancyhead[LE]{\footnotesize  Introduction}
\section{\xiaosan \textbf{Introduction and the Main results}}
In the present paper,  we study  the quasi-periodic Schr\"{o}dinger operator   $H=H_{\lambda v,\alpha,\theta}$ on  $   \ell ^2(\mathbb{Z})$:
\begin{equation}\label{G11}
    (H_{\lambda v,\alpha,\theta}u)_{n}=u_{n+1}+u _{n-1}+ \lambda v(\theta+n\alpha)u_{n},
\end{equation}
where $v: \mathbb{T}=\mathbb{R} / \mathbb{Z}\rightarrow \mathbb{R} $ is   the potential,  $\lambda$ is the coupling, $\alpha  $ is the frequency, and $\theta $ is the phase. In particular,  the almost Mathieu operator (AMO) is given by (\ref{G11}) with $v(\theta)=2\cos(2\pi \theta)$,
denoted by $H_{\lambda,\alpha,\theta}$.
\par

For $\lambda=0$, it is easy to verify  that   Schr\"{o}dinger operator  ($ \ref{G11} $)
 has   purely absolutely continuous spectrum ($[-2,2]$) by Fourier transform. We expect   the  property (has   purely absolutely continuous spectrum)
 preserves  under sufficiently small perturbation, i.e., $\lambda$ is small. Usually there are two  smallness about $ \lambda $. One is perturbative, meaning that the smallness   $ \lambda $ depends not only on the potential $v$, but also on the frequency $\alpha$; the other  is  non-perturbative, meaning that the smallness
condition   only depends on the potential $v$, not on $\alpha$.
 \par
 It is well known  that $H_{\lambda v,\alpha,\theta}$  has purely absolutely continuous spectrum for $\alpha \in \mathbb{Q}$ and all $\lambda$.
  Thus, unless  stated otherwise,   we always assume  $\alpha\in \mathbb{R}\backslash \mathbb{Q}$  in the present paper.
  We also   assume $v$ is real analytic in a strip of real axis from now on.
  \par
   The following notions are essential in the study of equation  (\ref{G11}).
  \par
  We say $\alpha \in \mathbb{R}\backslash \mathbb{Q} $ satisfies a Diophantine condition $\text{DC}(\kappa,\tau)$ with $\kappa>0$ and $\tau>0$,
if
$$ ||k\alpha||_{\mathbb{R}/\mathbb{Z}}>\kappa |k|^{-\tau}  \text{ for   any } k\in \mathbb{Z}\setminus \{0\},$$
where $||x||_{\mathbb{R}/\mathbb{Z}}=\min_{\ell \in \mathbb{Z}}|x-\ell| $.
Let $\text{DC}=\cup_{\kappa>0,\tau>0}\text{DC}( \kappa,\tau)$. We say $\alpha$  satisfies   Diophantine condition, if $\alpha\in \text{DC}$.
 \par
Let
 \begin{equation}\label{G12}
  \beta= \beta(\alpha)=\limsup_{n\rightarrow\infty}\frac{\ln q_{n+1}}{q_n},
 \end{equation}
 where $ \frac{p_n}{q_n} $ is  the continued fraction approximants   to $\alpha$.
   One usually calls set $\{\alpha \in \mathbb{R}\backslash \mathbb{Q}|\;\beta(\alpha)>0\}$    exponential regime and set $\{\alpha\in \mathbb{R}\backslash \mathbb{Q}|\;\beta(\alpha)=0\}$   sub-exponential regime. Notice that  $\beta(\alpha)=0$ for $\alpha\in DC$.
   \par
In $ \cite{E} $, Eliasson  treats  ($ \ref{G11} $)  as a dynamical systems problem--reducibility of associated cocycles. He
 shows that such cocycles are     reducible  for a.e. spectrum, and
gives good estimates for the non-reducible ones via a sophisticated  KAM-type methods, which breaks the limitations of the earlier KAM methods, for instance, the work of Dinaburg and Sinai \cite{Din}(they need exclude some parts of the spectrum).  As a result, Eliasson proves that   $H=H_{\lambda v,\alpha,\theta}$
has   purely absolutely continuous spectrum for all $\theta$,
 if  $\alpha\in DC$   and $|\lambda|<\lambda_0(\alpha,v)$\footnote{$\lambda_0(\ast)$ means $\lambda_0$   depends on $\ast$.  }. Clearly,  Eliasson's result is perturbative.
     \par
Bourgain and  Jitomirskaya established the  measure-theoretic version in non-perturbative regime, more precisely, they  proves that for a.e. $\alpha $ and  $\theta$, $H=H_{\lambda v,\alpha,\theta}$ ($H_{\lambda,\alpha,\theta}$) has  purely absolutely continuous spectrum
 if $|\lambda|<\lambda_0(v)$ ($\lambda<1$), see \cite{B2},\cite{BJ2},\cite{J} for some details.  They  approach   this by classical Aubry dulity  and the sharp  estimates of Green function
 in the regime of positive Lyapunov exponent.
  Bourgain    list a example  which suggests that the non-perturbative results in  multifrequency\footnote{ The
  quasi-periodic Schr\"{o}dinger operator in multifrequency($k$ dimension, $k\geq 2$)is given by  $( H_{\lambda v,\alpha,\theta}u)_n=u_{n+1} +u _{n-1}+ \lambda v(\theta+n\alpha)u_n,$  where $ v: \mathbb{T}^k= \mathbb{R}^k/ \mathbb{Z}^k\rightarrow \mathbb{R} $  is   the potential and
     $\alpha=(\alpha_1,\alpha_2,\cdots, \alpha_k)$  is such that
 $ 1,\alpha_1,\cdots, \alpha_k$ are independent over the rational numbers.} is wrong {\cite{B2}}.
 \par
 In $ \cite{AJ2} $,  Avila and Jitomirskaya   firstly  develop a quantitative version of  Aubry duality (Lemma \ref{Le33}) for  $\alpha\in DC$. As an application, they show that
  operator ($\ref{G11}$) has   purely absolutely continuous spectrum
in non-perturbative regime  for  $\alpha\in DC$ and all $\theta$, by reducing   non-perturbative regime to   Eliasson's perturbative regime.
In addition the sharp estimates of rotation number and      transfer matrix  (\cite{JL1},\cite{JL2}),
  Avila   prove that $ H_{\lambda v,\alpha,\theta}$ has
purely absolutely continuous spectrum in non-perturbative regime if $\beta(\alpha)=0$ \cite{A2}.
\par
 The present authors obtain the sharp estimate of rotation number in \cite{LIU1},  and extend  the quantitative version of  Aubry duality to  exponential
 regime\cite{LIU2}. Combining with Avila's arguments in \cite{A2}, we obtain the main  theorem in the present paper.
 \begin{theorem}\label{Th11}
 For irrational number $ \alpha$ such that $0< \beta(\alpha)<\infty$, if $v$ is  analytic in strip $|\Im x|<C\beta$, where $C$ is a large absolute
 constant, then
   there exists $\lambda_0=\lambda_0(v,\beta)$ such that
   $ H_{\lambda v,\alpha,\theta}$   has  purely absolutely continuous spectrum
 if $|\lambda|<\lambda_0$.
 \end{theorem}
\section{\xiaosan \textbf{ Preliminaries    }}
\subsection{\xiaosan \textbf{   Cocycles }}
  Denote by $ \text{SL}(2,\mathbb{C})$ the all complex  $2\times 2 $-matrixes  with determinant 1.
We say a function
  $f\in C^{\omega}(\mathbb{R}/\mathbb{Z},  \mathbb{C})$ if $f$ is well  defined in $ \mathbb{R}/\mathbb{Z}$, i.e., $f(x+1)=f(x)$,
  and $f$ is   analytic in
    a strip of  real axis. The definitions of $ \text{SL}(2,\mathbb{R})$ and  $  C^{\omega}(\mathbb{R}/\mathbb{Z},  \mathbb{R})$
    are similar to those of $ \text{SL}(2,\mathbb{C})$ and $C^{\omega}(\mathbb{R}/\mathbb{Z},\mathbb{C})$  respectively, except that the involved matrixes are real and the functions are real   analytic.
    \par
A $C^{\omega}$-cocycle in $ \text{SL}(2,\mathbb{C})$  is a pair $(\alpha,A)\in \mathbb{R} \times  C ^{\omega} (\mathbb{R}/\mathbb{Z}, \text{SL}(2,\mathbb{C})) $, where
$A  \in    C ^{\omega} (\mathbb{R}/\mathbb{Z}, \text{SL}(2,\mathbb{C}))$ means
      $A(x)\in \text{SL}(2,\mathbb{C})$      and the elements of $A$ are in $  C^{\omega}(\mathbb{R}/\mathbb{Z},  \mathbb{C})$.
       Sometimes, we say
     $ A$ a $C^{\omega}$-cocycle for short, if there is no ambiguity. Note that all functions, cocycles in the present paper are analytic in a strip of  real axis.
     Thus we often   do not mention the analyticity, for instance, we say $A$  a cocycle instead of   $C^{\omega}$-cocycle.
     \par
     Given two cocycles $(\alpha,A)$  and $(\alpha,A^{\prime})$, a conjugacy between them is a   cocycle
$ B \in C ^{\omega}(\mathbb{R}/\mathbb{Z},  \text{SL}(2,\mathbb{C}))$ such that
 \begin{equation}\label{G21}
   B(x+\alpha)^{-1}A(x)B(x)=A^{\prime}.
 \end{equation}
  The notion
of real conjugacy (between real cocycles) is the same as before, except that we ask for
 $ B \in C ^{\omega}(\mathbb{R}/\mathbb{Z}, \text{PSL}(2,\mathbb{R}))$, i.e., $B(x+1)= \pm B(x)$ and $\det B=1$.
We say that cocycle $(\alpha,A)$  is reducible if it is  conjugate to a constant cocycle.
\par
The Lyapunov exponent for the cocycle $A$ is given  by
 \begin{equation}\label{G22}
    L(\alpha,A)=\lim_{n\rightarrow\infty} \frac{1}{n}\int_{\mathbb{R}/\mathbb{Z}} \ln \| A_n(x)\|dx,
 \end{equation}
 where
  \begin{equation}\label{G23}
     A_n(x) = A(x+(n-1)\alpha)A(x+(n-2)\alpha)\cdots A(x).
\end{equation}
We say cocycle $(\alpha,A)$ is bounded if $\sup_{n\geq0,x\in \mathbb{R}}||A_n(x)|| <\infty$.
 \par
      We now consider the  quasi-periodic Schr\"{o}dinger operator  \{$ H_{\lambda v,\alpha,\theta} \}_{\theta\in \mathbb{R}}$. It is easy
       to verify that the spectrum of $ H_{\lambda v,\alpha,\theta} $
  does not depend on $ \theta$ for $\alpha\in  \mathbb{R}\backslash\mathbb{Q}$, thus we denote by $  \Sigma_{\lambda v,\alpha}$.
  \par
 Let
 $$
 S_{\lambda v,E}=
 \left(
     \begin{array}{cc}
       E- \lambda v & -1 \\
       1 & 0 \\
     \end{array}
   \right).
   $$
We call $(\alpha,S_{\lambda  v,E})$  Schr\"{o}dinger cocycle.
\par
Fix Schr\"{o}dinger operator $  H_{\lambda v,\alpha,\theta}   $, we define the Aubry dual model by $\hat{H}= \hat{H}_{\lambda v,\alpha,\theta}  $,
\begin{equation}\label{G24}
  (\hat{H}\hat{u}) _ n = \sum _{k\in \mathbb{Z}}\lambda \hat{v}_ k \hat{u}_{n-k}+2\cos(2\pi \theta+n\alpha)\hat{u}_n,
\end{equation}
where $\hat{v}_ k$ is the Fourier coefficients of potential $v$.
If $\alpha\in\mathbb{R} \backslash  \mathbb{Q}$,  the spectrum of $ \hat{H}_{\lambda v,\alpha,\theta} $   is also $\Sigma_{\lambda v ,\alpha}$   \cite{GJLS}.
  \subsection{\xiaosan \textbf{The rotation number}}
 Let
$A(\theta)=\left(
             \begin{array}{cc}
               a(\theta) & b(\theta) \\
               c(\theta)& d (\theta)\\
             \end{array}
           \right)
$, we define  the map $T_{\alpha,A}:(\theta,\varphi)\in \mathbb{T}\times \frac{1}{2}\mathbb{T} \mapsto
(\theta+\alpha,\varphi_{\alpha,A}(\theta,\varphi))\in \mathbb{T}\times \frac{1}{2}\mathbb{T},$
 with $\varphi_{\alpha,A}= \frac{1}{2\pi}\arctan (\frac{c(\theta)+d(\theta)\tan2\pi\varphi}{a(\theta)+b(\theta)\tan2\pi\varphi})$,
 where $\mathbb{T}=\mathbb{R}/\mathbb{ Z}$.
Assume now that     A :$ \mathbb{R}/   \mathbb{ Z}\rightarrow \text{SL}(2,\mathbb{R})$ is homotopic to the identity,
then  $T_{\alpha,A} $ admits a continuous lift $\tilde{T}_{\alpha,A}:(\theta,\varphi)\in \mathbb{R}\times  \mathbb{R} \mapsto
(\theta+\alpha,\tilde{\varphi}_{\alpha,A}(\theta,\varphi))\in \mathbb{R}\times  \mathbb{R} $  such that
$\tilde{\varphi}_{\alpha,A}(\theta,\varphi)  \mod  \frac{1}{2}\mathbb{Z}= \varphi _{\alpha,A}(\theta,\varphi  ) $
and $\tilde{\varphi}_{\alpha,A}(\theta,\varphi) -\varphi$ is well defined on $\mathbb{T}\times \frac{1}{2}\mathbb{T}$.
 The number $\rho(\alpha,A)=\limsup_{n\rightarrow \infty}\frac{1}{n}(p_2\circ \tilde{T}^n_{\alpha,A}(\theta,\varphi)-\varphi)\mod \frac{1}{2}\mathbb{Z},$
 does not depend on the choices of $ \theta$ and $\varphi$,
where $p_2(\theta,\varphi)=\varphi,$ and is called the   rotation number
of $(\alpha,A)$ $ {\cite{Her} ,\cite{JM}}$.
\par
It's easy to see that  Schr\"{o}dinger   cocycle is  homotopic to the identity,
  and let
$\rho_{\lambda v,\alpha}(E)\in[0,\frac{1}{2}]$ be the rotation number of Schr\"{o}dinger cocycle $(\alpha,S_{\lambda v,E})$.
\subsection{Spectral measure and the    integrated density of states}
Let $H$ be a bounded self-adjoint operator on $\ell ^2(\mathbb{Z})$. Then
  $(H-z)^{-1}$ is analytic
in  $ \mathbb{C} \backslash \Sigma( H)$, where $\Sigma( H)$ is the spectrum of $H$,
and we have for $ f\in \ell^2$
\begin{equation*}
  \Im\langle(H-z)^{-1}f,f\rangle=\Im z \cdot ||(H-z)^{-1}f||^2,
\end{equation*}
where $\langle\cdot,\cdot\rangle$ is the usual inner product in $\ell ^2(\mathbb{Z})$.
Thus
\begin{equation*}
 \phi_f(z)=   \langle(H-z)^{-1}f,f\rangle
\end{equation*}
is an analytic function in the upper half plane with $\Im \phi_f\geq 0$ ( $\phi_f$ is a so-called Herglotz function).
\par
Therefore one has a representation
\begin{equation}\label{G25}
 \phi_f(z)=   \langle(H-z)^{-1}f,f\rangle=\int_{\mathbb{R}} \frac{1}{x-z}d\mu^{f}(x),
\end{equation}
where $\mu^{f} $ is the spectral measure   associated to vector $f$.
\par
Denote by $\mu^{f}_{\lambda v ,\alpha,\theta} $
the spectral measure of operator  $  H_{\lambda  v,\alpha,\theta}   $ and vector $f$ as before.
The  integrated density of states   $N_{\lambda v  ,\alpha}$ is obtained by   averaging the spectral measure  $\mu_{\lambda v , \alpha,\theta}^{e_0}$ with respect to $\theta$, where
$e_0$ is the Dirac mass at $0 \in \mathbb{Z}$, i.e.,
\begin{equation} \label{G26}
  N_{\lambda v,\alpha}(E)=\int_{\mathbb{R}/\mathbb{Z}} \mu^{e_0}_{\lambda v  ,\alpha,\theta}(-\infty,E]d\theta.
\end{equation}
\par
Between the integrated density  of   states $N_{\lambda v,\alpha}(E)$ and the rotation number  $\rho_{\lambda v,\alpha}(E)$,
there is   the following relation $ {\cite{Jo1}}$:
 \begin{equation}\label{G27}
    N_{\lambda v,\alpha}(E)=1-2\rho_{\lambda v,\alpha}(E).
 \end{equation}

\section{ Some known results }
Let $ \mu_{\lambda v,\alpha,\theta}=\mu_{\lambda v,\alpha,\theta}^{e_{-1}}+\mu_{\lambda v,\alpha,\theta}^{e_0}$,
where
$e_i$ is the Dirac mass at $i \in \mathbb{Z}$. For simplicity, sometimes we  drop the parameters dependence, for example,
   replacing $\mu_{\lambda v,\alpha,\theta}$    with $\mu$.
Fix  $A=
 S_{\lambda v,E}=
 \left(
     \begin{array}{cc}
       E- \lambda v & -1 \\
       1 & 0 \\
     \end{array}
   \right).
   $    Below,  $C$ is a large absolute  constant and $c$ is a small
 absolute  constant, which may change through the arguments, even when appear in the same formula.  Denote by
   $C_{\star}$  $(c_{\star})$ a large(small) constant depending on $\lambda, v,\alpha$. Let  $||\cdot||$  be the Euclidean norms, and denote
 $||f||_\eta=\sup _{|\Im x|<\eta} ||f(x)| |$,  $||f||_0=\sup _{x\in \mathbb{R}} ||f(x) ||$.
 \begin{lemma}$(\text{ Lemma } 2.4,\cite{A2} )$\label{Le31}
  Let $\mathfrak{B}$ be the set of $E\in \mathbb{R}$ such that  the cocycle $(\alpha,A)$ is bounded,  then
  $\mu| _{ \mathfrak{B}}$ is absolutely continuous.
 \end{lemma}
\begin{lemma}$(\text{ Lemma } 2.5,\cite{A2} )$\label{Le32}
  We have $\mu(E-\varepsilon,E+\varepsilon)\leq C \varepsilon\sup_{0\leq s\leq C\varepsilon^{-1}}||A_s||_0^2$.
 \end{lemma}
Given    $\epsilon_0>0$,  we say $k$ is an $\epsilon_0$-resonance for $\theta$, if
$ \| 2\theta-k\alpha\|_{\mathbb{R}/\mathbb{Z}}\leq e^{-\epsilon_0|k|}$ and
$\| 2\theta-k\alpha\|_{\mathbb{R}/\mathbb{Z}}=\min_{|j|\leq|k|} \| 2\theta-j\alpha\|_{\mathbb{R}/\mathbb{Z}}$.
\par
Clearly, $0\in \mathbb{Z}$ is  an $\epsilon_0$-resonance.
We order the $\epsilon_0$-resonances $0=|n_0|<|n_1|\leq|n_2|\cdots$. We say  $\theta$ is $\epsilon_0$-resonant if the
set of $\epsilon_0$-resonances is infinite.
\begin{lemma}$(\text{Theorem }3.3, \cite{AJ2})$\label{Le33}
 If $E\in \Sigma_{\lambda v,\alpha} $, then there exists $ \theta\in \mathbb{R}$ and a bounded solution  of
 $\hat{H}_{\lambda v, \alpha,\theta}\hat{u}=E\hat{u}$ with
 $\hat{u}_ 0  =1$ and $|\hat{u}_ k |\leq 1$.
\end{lemma}
Fix $\alpha\in \mathbb{R}\backslash \mathbb{Q}$ such that $0<\beta(\alpha) <\infty$. Let $\epsilon_0=C_1^2\beta$, where
$C_1$ is a large absolute constant, which is much larger than any absolute constant $C$, $c^{-1}$ emerging in the present paper.  Set $ h_1=C_1 \beta$, $ h_2=C_1^3\beta$. Fix $E\in \Sigma_{\lambda v,\alpha}$ below,  and
choose some    $\theta = \theta(E)$     given by Lemma $ \ref{Le33}$.  Denote $\{n_j\} $ all the  $\epsilon_0$-resonances for $\theta(E)$.

By the present authors's arguments in \cite{LIU1},\cite{LIU2},
if $v$ is  analytic in  strip $|\Im x|<C_2\beta$, where $C_2$ is a large absolute
 constant,
then there exists $ \lambda_0=\lambda_0(v,\beta)>0$  such that the following theorems hold for $|\lambda|<\lambda_0$.
\begin{theorem}$(\text{Lemma  }4.3, \cite{LIU2})$\label{Th34}
 The  Lyapunov exponent vanishes on $\Sigma_{\lambda v,\alpha}$, i.e., $L(\alpha,S_{\lambda v,E})=0$ for all $E\in \Sigma_{\lambda v,\alpha}$.
\end{theorem}
\begin{theorem} $(\text{Theorem } 5.6, \cite{LIU2})$ \label{Th35}
We have the following estimate,
\begin{equation}\label{G32}
    ||A_s||_0\leq C_{\star} e^{C \beta n},0\leq s\leq  c_{\star}e^{c\epsilon_0 n}.
\end{equation}
 \begin{theorem} $(\text{Corollary } 6.2, \cite{LIU2})$\label{Th36}
  The integrated density of states  of $H_{\lambda v,\alpha,\theta}$ is 1/2-H\"{o}lder continuous, that is $N_{\lambda v,\alpha}(J)\leq C_{\star}|J|^{1/2}$ for any interval $J\subset \mathbb{R}$.
  \end{theorem}
\end{theorem}
\begin{theorem} $(\text{Theorem } 4.14, \cite{LIU1})$\label{Th37}
 If $\theta=\theta(E)$  has a   $\epsilon_0$-resonance $ n_j  $,
then there exists  $m_j$ with $|m_j|\leq C |n_j| $  such that $ || 2\rho_{\lambda v,\alpha}(E) -m_j\alpha ||_{\mathbb{R}/\mathbb{Z}}\leq C_{\star}e^{-\epsilon_0 |n_j|}$, or equivalently $ || N_{\lambda v,\alpha}(E) -m_j\alpha ||_{\mathbb{R}/\mathbb{Z}}\leq C_{\star}e^{-\epsilon_0 |n_j|}$ by (\ref{G27}).
  \end{theorem}
  \begin{theorem}  $(\text{ Theorem  } 5.8, \cite{LIU1})$ \label{Th38}
     If $\theta(E)$ is not  $ \epsilon_0$-resonant, then cocycle $A$ is reducible.
     \begin{remark}
  In \cite{LIU1}, the present authors only prove Theorem  \ref{Th37}, \ref{Th38} for AMO  by   quantitative version of  Aubry duality
 in exponential
 regime. For the general quasi-periodic Schr\"{o}dinger operator, the proof is
  similar if we use the  quantitative version of  Aubry duality for general potential $v$ in \cite{LIU2}.
  \end{remark}
  \end{theorem}

\section{Proof of   Theorem \ref{Th11}}
  \begin{lemma}\label{Le41}  For $0<\varepsilon<1$,  $N_{\lambda v,\alpha}(E+\varepsilon) -N_{\lambda v,\alpha}(E-\varepsilon)\geq c_{\star}\varepsilon ^{2}$.
 \end{lemma}
 \textbf{Proof:} The lemma can be proved directly   by   Theorem \ref{Th34} and \ref{Th36}. See the proof of  Lemma 3.11 in $ \cite{A2}$ for details.
 \par
 \textbf{ Proof of theorem  } $\ref{Th11}$:
  It is well known that it suffices to prove that $\mu$ is   absolutely continuous.
   Let $\mathfrak{B}$ be given by Lemma \ref{Le31}. Thus  it   suffices to show $\mu(\Sigma_{\lambda v,\alpha}\backslash \mathfrak{B} )=0$.
    Let $\mathscr{R} $   be the set of $E\in  E_{\lambda v,\alpha}$  such that  $A$ is reducible.
  We have $\mu(  \mathscr{R}\backslash\mathfrak{B} )=0$, since $ \mathscr{R}\backslash\mathfrak{B} $ is a countable set and there is no eigenvalue in $ \mathscr{R} $(see p.16 in \cite{A2} for details).
  Thus to prove the Theorem \ref{Th11}, it is sufficient to show that $\mu( \Sigma_{\lambda v,\alpha}\backslash \mathscr{R})=0$.
 \par
Let $K_m\subset \Sigma _{\lambda v,\alpha}$, $m\geq1$  be the set of $E$ such that there exists $\theta (E) \in \mathbb{R} $  given by Lemma $ \ref{Le33}$  with a resonance $ 2^m\leq|n_j|\leq 2^{m+1}$. We will show that  $\sum \mu(\overline{K_m})<\infty$. By Theorem $ \ref{Th38}$ $  \Sigma_{\lambda v,\alpha}\backslash \mathscr{R}\subset \limsup K_m$, then
 $\mu( \Sigma_{\lambda v,\alpha}\backslash \mathscr{R})=0$  by the fact
 $\sum \mu(\overline{K_m})<\infty$   and the Borel-Cantelli Lemma.
\par
 For every  $E\in K_m$, let $J_m(E) $ be an open $\epsilon_m=C_{\star}e^{-c \epsilon_0 2^{m-1}}$
neighborhood of $E$.  By   $(\ref{G32})$,
\begin{equation}\label{G41}
 \sup_{0\leq s\leq C\epsilon_m^{-1}}||A_s||_0\leq C_{\star} e^{C\beta2^m} .
\end{equation}
  Take a finite subcover $\overline{K}_m \subset  \cup_{j=0}^r J_m(E_j)$.
Refining this subcover if necessary, we may assume that every $x\in\mathbb{R}$ is contained
in at most 2 different $J_m(E_j)$.
\par
 By lemma $ \ref{Le41}$, $ N(J_m(E))\geq c_{\star}|J_m(E)|^2\geq C_{\star}e^{- c\epsilon_0 2^m}$. By Theorem $\ref{Th37}$, if $E \in K_m$ then
  $||N(E)-k\alpha||_{\mathbb{R}/\mathbb{Z}}\leq C_{\star} e^{-\epsilon_0 2^m} $ for some
   $|k|\leq C2^m$, so there are at most $C_{\star}2^m$ intervals $J_m(E_j)$, i.e., $r\leq C_{\star}2^m$. Thus by (\ref{G41}) and Lemma \ref{Le32},
   \begin{equation}\label{G42}
    \mu(\overline{K}_m)\leq \sum_{j=0}^r\mu (J_m(E_j))\leq C_{\star} 2^m e^{C\beta 2^m}e^{-c\epsilon_0 2^{m-1}},
   \end{equation}
which implies  $\sum_m \mu(\overline{K_m})<\infty$.$\qed$

\par
Next, we will prove that the integrated density of states is absolutely continuous in perturbative regime for all $\alpha$ satisfying
$ 0<\beta(\alpha)<\infty$.  We need   a lemma  first.
\begin{lemma}$(\text{Corollary  }1, \cite{Dam})$ \label{Le42}
If the Lyapunov exponent vanishes on $\Sigma_{\lambda v,\alpha}$, then $H_{\lambda v,\alpha,\theta}$ has purely absolutely continuous spectrum
for almost $\theta$  if and only if the integrated density of states $N_{\lambda v,\alpha}(E)$ is absolutely continuous.
\end{lemma}

 \begin{theorem}\label{Th43}
 For irrational number $ \alpha$ such that $0< \beta(\alpha)<\infty$, if $v$ is  analytic in strip $|\Im x|<C_2\beta$, where $C_2$ is a large absolute
 constant, then
   there exists $\lambda_0=\lambda_0(v,\beta)$ such that
   the integrated density of states $N_{\lambda v,\alpha}(E)$ is absolutely continuous if $|\lambda|< \lambda_0$.
 \end{theorem}
\textbf{Proof:} Using Theorem \ref{Th34} and Lemma \ref{Le42}, $N_{\lambda v,\alpha}(E)$ is absolutely continuous  if and only if $H_{\lambda v,\alpha,\theta}$ has purely absolutely continuous spectrum
for almost every $\theta$.   Together  with Theorem \ref{Th11}, we finish the proof.

           \begin{center}
           
             \end{center}
  \end{document}